\documentstyle[12pt,amsfonts,amssymb,leqno]{article}

\newfont\got{eufm10}

\newtheorem{proposition}{Proposition}[section]
\newtheorem{thm}[proposition]{Theorem}

\newtheorem{lemma}[proposition]{Lemma}
\newtheorem{defn}[proposition]{Definition}

\newtheorem{convention}[proposition]{Convention}
\newtheorem{claim}[proposition]{Claim}
\renewcommand{\thefootnote}{\alph{footnote}}

\newcounter{secnum}

\begin{document}
\setcounter{section}{+0}

\begin{center}
{\Large \bf Homogeneously Souslin sets in small inner models}
\end{center}

\begin{center}
\renewcommand{\thefootnote}{\fnsymbol{footnote}}
{\large Peter Koepke}${}^{a}$
\renewcommand{\thefootnote}{arabic{footnote}}
\renewcommand{\thefootnote}{\fnsymbol{footnote}}
{\large and Ralf Schindler}${}^{b}$
\end{center}
\begin{center} 
{\footnotesize
${}^{a}${\it 
Mathematisches Institut,  
Universit\"at Bonn, 
Beringstra\ss e 1, 
Bonn, Germany 
}\\
${}^b${\it Institut f\"ur Formale Logik, Universit\"at Wien, 1090 Wien, Austria}}
\end{center}

\begin{center}
{\tt koepke@math.uni-bonn.de, 
rds@logic.univie.ac.at}\\
\end{center}

\begin{abstract}
\noindent
We prove that every homogeneously Souslin set is coanalytic
provided that either (a) $0^{\rm long}$ does not exists, or else
(b) $V = K$, where $K$ is the core model below a $\mu$-measurable cardinal. 
\end{abstract}

\section{Homogeneously Souslin sets.}

In this paper we shall deal with homogeneously Souslin sets of reals,
or rather with sets of reals which admit an $\omega$-closed embedding
normal form.

\begin{defn} (Cf.~\cite[p.~92]{tony-john}.)
Let $A \subset {}^\omega\omega$. Let $\alpha \in {\rm OR}$.
We say that $A$ has an $\alpha$-{\em closed 
embedding normal form} if and only if the following holds true.
There is a commutative
system $$((M_s \colon s \in 
{}^{<\omega}\omega),(\pi_{s t} \colon s, t 
\in {}^{<\omega}\omega, s \subset t))$$
such that $M_0 = V$, each $M_s$ is an inner model of ${\sf ZFC}$ with ${}^\alpha M_s
\subset M_s$, each
$\pi_{s t} \colon M_s \rightarrow M_t$ is an elementary embedding, 
and if $x \in {}^\omega \omega$ and 
$(M_x , (\pi_{x \upharpoonright n , x} \colon
n<\omega))$ is the direct limit of $((M_{x \upharpoonright n} \colon n<\omega) , 
(\pi_{x \upharpoonright n , x \upharpoonright m} \colon
n \leq m <\omega))$ then $$x \in A \ \Leftrightarrow M_x {\rm \ is \ wellfounded}.$$ 
\end{defn}

As we shall not need it here, we do not repeat the definition of 
the concept of being homogeneously Souslin in this paper (cf.~\cite[p.~87]{tony-john}). We just remind the
reader of the following facts.

\begin{lemma}\label{fact} 
Let $A \subset {}^\omega\omega$. 

(1) If $A$ is coanalytic and if $\kappa$ is a measurable cardinal
then $A$ is $\kappa$-homogeneously Souslin
(cf.~\cite{tony}, \cite[Theorem 2.2]{tony-john}).

(2) If $A$ is
$\kappa$-homogeneously Souslin, where
$\kappa$ is a (measurable) cardinal, then $A$ 
is determined (cf.~\cite[Theorem 2.3]{tony-john}) and
has a $\kappa$-closed embedding
normal form (cf.~\cite[p.~92]{tony-john}).

(3) If $A$ has a $2^{\aleph_0}$-closed embedding normal form
then $A$ is homogeneously Souslin (cf.~\cite[Lemma 2.5]{katrin}, 
\cite[Theorem 5.2]{peter}).
\end{lemma}

Our aim is to prove a converse to Lemma \ref{fact} (1) 
under appropriate anti-large cardinal hypotheses.

%
%
%
\begin{defn} 
A cardinal $\kappa$ is called $\mu${\em -measurable} if there is an
embedding $\pi \colon V \rightarrow M$ 
such that $M$ is transitive,
$\kappa = {\rm crit}(\pi)$, and 
$\{ X \subset \kappa | \kappa \in \pi(X) \} \in M$ (cf.~\cite{bill}).

We say that $0^\P$ does not exist if for every iterable premouse
${\cal M}$, if $E^{\cal M}_\nu \not= \emptyset$ then 
${\cal M}||{\rm crit}(E^{\cal M}_\nu)$ is a model of ``there is no strong 
cardinal (as being witnessed by the extenders from the
${\cal M}$-sequence)'' (cf.~\cite{martin}).

Suppose that $0^\P$ does not exist, and let $K$ denote the core model
(cf.~\cite[Chap.~8]{martin}). We say that $K$ is {\em below a $\mu$-measurable
cardinal} if $K \models$ ``there is no $\mu$-measurable cardinal.''
\end{defn}

If $K$ is below a $\mu$-measurable cardinal then every total extender
on the $K$-sequence has exactly one generator.

\begin{defn} 
We say that $0^{\rm long}$ does not exist if for every iterable premouse
${\cal M}$, if we let $A$ be the set of critical points of the total measures
from the ${\cal M}$-sequence then $A = \emptyset$ or else
${\rm otp}(A) < {\rm min}(A)$
(cf.~\cite{peter}).
\end{defn}

We can now state the main results of our paper.  

\begin{thm}\label{thm1} Suppose that $0^\P$ does not exist, and $K$ is
below a $\mu$-measurable cardinal.
Suppose that $V=K$. Let $A \subset {}^\omega\omega$ have an $\omega$-closed 
embedding normal form. Then $A$ is coanalytic.
\end{thm}

\begin{thm}\label{thm2}
Suppose that $0^{\rm long}$ does not exist.
Let $A \subset {}^\omega\omega$ have an $\omega$-closed 
embedding normal form. Then $A$ is coanalytic.
\end{thm}

Our main technical tool will be the concept of a ``shift map.''
Shift maps will be defined in the next section, where
we shall also show that if $K$ is below a $\mu$-measurable cardinal
then any elementary embedding from one universal weasel into another
one is a shift map. The final section will prove Theorems \ref{thm1}
and \ref{thm2}.

As to prerequisites, we shall assume familiarity with 
the core model theory as presented in \cite[Chap.~8]{martin}.

\section{Shift maps.}
In this section we shall prove a result on shift maps. This result is more 
general than what we would need in order to prove the main theorems, 
but it might be interesting in its own right. 

\begin{convention}
Let $\zeta \leq \theta$, and let 
$\varphi \colon \zeta+1 \rightarrow \theta+1$ be strictly
monotone, i.e., $\varphi(\alpha') > \varphi(\alpha)$ for $\alpha' > \alpha$. We then
let ${\varphi}^-$ denote the partial map from $\theta+1$ to $\zeta+1$
which is defined by 
${\varphi}^-(\beta) =$
the least $\alpha$ such that $\varphi(\alpha) \geq \beta$.
\end{convention}

Notice that in this situation, $\varphi^-$ is total (i.e., 
${\rm dom}(\varphi^-) = \zeta+1$) if and only if $\varphi(\zeta) =
\theta$.

\begin{defn}\label{mainthm}
Let $\pi \colon W \rightarrow W'$ be an elementary embedding, where both
$W$ and $W'$ are weasels.
We say that $\pi$ is a {\em shift map} provided the following holds
true.

There is a weasel $W_0$, 
there are non-dropping iterations
${\cal T}$ and ${\cal U}$ of $W_0$ 
with ${\rm lh}({\cal T}) = \zeta+1$ and ${\rm lh}({\cal U}) = \theta+1$, where 
$\zeta \leq \theta \leq {\rm OR}$, there is a strictly monotone 
$\varphi \colon \zeta+1 \rightarrow \theta+1$
with $\varphi(\zeta) = \theta$, and there are elementary embeddings
$\pi_\beta \colon {\cal M}_{\varphi^-(\beta)}^{\cal T} \rightarrow {\cal
M}_\beta^{\cal U}$ for $\beta \leq \theta$
such that:

(a) $\pi = \pi_{\theta}$,

(b) $\pi_{\beta'} \circ \pi_{\varphi^-(\beta) \varphi^-(\beta')}^{\cal T} =
\pi_{\beta \beta'}^{\cal U} \circ \pi_\beta$ whenever $\beta \leq \beta' \leq
\theta$, 

(c) $\pi_{\varphi(\alpha)+1} \upharpoonright {\rm lh}(E_\alpha^{\cal T}) =
\pi_{\varphi(\alpha)} \upharpoonright {\rm lh}(E_\alpha^{\cal T})$
whenever $\alpha + 1 \leq \zeta$, and

(d) $E_{\varphi(\alpha)}^{\cal U} = \pi_{\varphi(\alpha)}(E_\alpha^{\cal T})$
whenever $\alpha + 1 \leq \zeta$. 
\end{defn}

We shall prove:

\begin{thm}\label{main-theorem} 
Suppose that $0^\P$ does not exist, and $K$ is
below a $\mu$-measurable cardinal.
Let $\pi \colon W \rightarrow W'$ be an elementary embedding, where both
$W$ and $W'$ are universal weasels. Then $\pi$ is a shift map. 
\end{thm}

Let us state some terminology, assuming that
$0^\P$ does not exist, 
before commencing with the proof of Theorem 
\ref{main-theorem}. Let $W$ be a weasel and let $\Gamma \subset {\rm OR}$ be thick in
$W$ (cf.~\cite[p.~214f.]{martin}). Let $\alpha$ be an ordinal. We
shall write $H^W_\Gamma(\alpha)$ for the set of all $\tau^W({\vec
\gamma},{\vec \epsilon})$, where $\tau$ is a Skolem term, ${\vec \gamma} \in
[\alpha]^{<\omega}$, and ${\vec \epsilon} \in [\Gamma]^{<\omega}$.
We shall write $H^W(\alpha)$ for the intersection of all
$H_\Gamma^W(\alpha)$ where $\Gamma$ is thick in $W$.
$W$ has the {\em definability property at} $\alpha$ just in case that 
$\alpha \in H^W(\alpha)$ (cf.~\cite[Definition 4.4]{CMIP}).

Let ${\cal M}$ be a premouse, and let $\alpha \leq {\cal M} \cap {\rm OR}$. We say
that ${\cal M}$ is $\alpha$-{\em very sound} 
if there is a weasel $W$ such that ${\cal M}
\triangleleft W$ and ${\cal M} \subset H^W(\alpha)$.
In this case, $W$ is called an 
$\alpha$-{\em very soundness witness} for ${\cal M}$.

A premouse ${\cal M}$ is {\em strong} if there is a universal
weasel $W$ with ${\cal M} \triangleleft W$ (cf.~\cite[p.~212]{martin}).

Let ${\cal M}$ be strong, and let $W \triangleright {\cal M}$ be a witness to
this. By \cite[Theorem 7.4.9 and p.~275]{martin} there is a
normal non-dropping iteration ${\cal T}^*$ of
$K$ such that $W = {\cal M}_\infty^{{\cal T}^*}$. Therefore, if we let
$({\cal U},{\cal T})$
denote the coiteration of ${\cal M}$ with $K$ then ${\cal U}$ is trivial and ${\cal
T}$ is simple.
The following is part of the folklore.

\begin{lemma}\label{folklore} Suppose that $0^\P$ does not exist, and $K$ is
below a $\mu$-measurable cardinal.
Let $W$ be a universal weasel, let $\kappa$ be a cardinal of $W$,
and let ${\cal M} = W||\kappa^{(+3)W}$.
Then ${\cal M}$ is $(\kappa+1)$-very sound. Moreover, 
${\cal M}$ is $\kappa$-very sound if and only if the following holds true:
if ${\cal T}$ is the normal non-dropping iteration of $K$ which arises from the
comparison with ${\cal M}$ then no $E_\alpha^{\cal
T}$ has critical point $\kappa$.
\end{lemma}

Lemma \ref{folklore} is no longer true if $K$ is not assumed to be
below a $\mu$-measurable cardinal (cf.~\cite[p.~29, Example 4.3]{CMIP}).

\begin{lemma}\label{folklore2} Suppose that $0^\P$ does not exist, and $K$ is
below a $\mu$-measurable cardinal.
Let both $W$ and $W'$ be universal weasels. 
Let $W||\mu = W'||\mu$, where $\mu$ is
either a limit cardinal or else a double successor 
cardinal in both $W$ and $W'$.
Then $W||\mu^{+W} = W'||\mu^{+W'}$.
\end{lemma}

The following will be used towards the end of the proof of
Theorem \ref{main-theorem}.

\begin{lemma}\label{observation} 
Suppose that $0^\P$ does not exist, and $K$ is
below a $\mu$-measurable cardinal.
Let $W$, $W'$ be weasels, and let $\pi \colon W \rightarrow W'$ be elementary. Let
$E^W_\alpha \not= \emptyset$ with $\kappa = {\rm crit}(E^W_\alpha)$. Let $\pi$ be
continuous at $\kappa^{+W}$, and let $E^{W'}_\beta \not= \emptyset$ be such that
$$\pi {\rm " } E^W_\alpha \subset E^{W'}_\beta.$$ Then $\beta = \pi(\alpha)$, i.e., 
$E^{W'}_\beta = \pi(E^W_\alpha)$.
\end{lemma}

{\sc Proof} of Lemma \ref{observation}.
Set $\lambda = \pi(\kappa)$. We must have $\lambda = {\rm
crit}(E_\beta^{W'})$. Let $X \in E_\beta^{W'}$. We aim to prove that $X \in
E_{\pi(\alpha)}^{W'}$.

As $\pi$ is
continuous at $\kappa^{+W}$ we may pick some $\gamma < \kappa^{+W}$ such that $X \in
W'||\pi(\gamma)$. Let $(X_i \colon i<\kappa) \in W$ be such that $E_\alpha^W \cap
W||\gamma = \{ X_i \colon i<\kappa \}$. Let $Y = \Delta_{i<\kappa} X_i$ be the
diagonal intersection. 
As $E_\alpha^W$ is normal, $Y \in E_\alpha^W$, and thus 
$\pi(Y) \in E_\beta^{W'}$ by $\pi {\rm " } E^W_\alpha \subset E^{W'}_\beta$.
Moreover, for each $Z \in {\cal P}(\kappa) \cap W||\gamma$ there is some $\xi <
\kappa$ such that $Y \setminus \xi \subset Z$ or $Y \setminus \xi \subset \kappa
\setminus Z$ (we'll have the former if and only if $Z \in
E_\alpha^W$). By elementarity, hence,
for each $Z \in {\cal P}(\lambda) \cap W'||\pi(\gamma)$ there is some $\xi <
\lambda$ such that $\pi(Y) \setminus \xi \subset Z$ or $\pi(Y) 
\setminus \xi \subset \lambda
\setminus Z$.

In particular, there is some $\xi <
\lambda$ such that $\pi(Y) \setminus \xi 
\subset X$ or $\pi(Y) \setminus \xi \subset \lambda \setminus X$.
As $\pi(Y) \in E_\beta^{W'}$ and $X \in E_\beta^{W'}$ we cannot have that
$\pi(Y) \setminus \xi \subset \lambda \setminus X$. Therefore, $\pi(Y) \setminus \xi 
\subset X$.
But because $\pi(Y) \in E_{\pi(\alpha)}^{W'}$ we then get
$X \in E_{\pi(\alpha)}^{W'}$, as desired. 

\hfill $\square$ (Lemma \ref{observation})

\bigskip
We are now ready to prove the key result of this section.

\bigskip
{\sc Proof} of Theorem 
\ref{main-theorem}. As $W$ and $W'$ are universal, 
there are non-dropping iterations
${\cal T}$ and ${\cal U}$ of $K$ with
${\rm lh}({\cal T}) = \zeta+1$ and ${\rm lh}({\cal U}) = \theta+1$ for some
$\zeta$, $\theta \leq {\rm OR}$ such that $W = {\cal M}^{\cal T}_\zeta$ and
$W' = {\cal M}^{\cal U}_\theta$ (cf.~\cite[Theorem 7.4.9 and p.~275]{martin}). 

\begin{claim}\label{claim1}
Let $\alpha+1 \leq \zeta$,  
$\kappa = {\rm crit}(E_\alpha^{\cal T})$, and $\lambda = \pi(\kappa)$. Then
$\lambda = {\rm crit}(E_\beta^{\cal U})$ for some $\beta$ with $\beta+1
\leq \theta$. 
\end{claim}

{\sc Proof} of Claim \ref{claim1}. Set $\nu = {\kappa}^{(+3)W}$. 
Let $Q$ and $\Gamma$ be such that
$Q \triangleright W||\nu$ is a universal weasel, $\Gamma$ is thick in
$Q$, and
$\kappa \notin H^Q_\Gamma(\kappa)$. 
In particular, $\Gamma$ witnesses 
that $Q$ does not have the definability property at $\kappa$.
Let $E = E_{\pi \upharpoonright (W||\nu)}$ be the long extender derived
from $\pi \upharpoonright (W||\nu)$, and let
$${\tilde \pi} \colon Q \rightarrow_E {\tilde Q} = Ult(Q;E)$$
be the ultrapower of $Q$ by $E$. 
Of course, $\pi \upharpoonright (W||\nu) = {\tilde \pi} 
\upharpoonright (W||\nu)$.
Because ${\tilde Q}$ is universal, ${\tilde \pi}(\nu) = \pi(\nu)$
by Lemma \ref{folklore2}.
Thus $W'||\pi(\nu) \triangleleft {\tilde Q}$.
Let $${\tilde \Gamma} = \Gamma \cap \{ \epsilon \colon {\tilde \pi}(\epsilon) =
\epsilon \}.$$ We'll have that ${\tilde \Gamma}$ is
thick in both $Q$ and ${\tilde Q}$.

Now suppose that there is no $\beta$ with $\beta+1
\leq \theta$ and $\lambda = {\rm crit}(E_\beta^{\cal U})$.
Then ${\tilde Q}$ has the definability property at $\lambda$, and 
there is a Skolem term $\tau$ and there are
${\vec \gamma} \in \lambda$ and 
${\vec \epsilon} \in {\tilde \Gamma}$
such that $\lambda = \tau^{\tilde Q}({\vec \gamma},{\vec \epsilon})$. That is,
$${\tilde Q} \models \exists {\vec {\bar \gamma}} < \lambda \ \lambda = 
\tau({\vec {\bar \gamma}},{\vec \epsilon}).$$
Therefore, by using the map ${\tilde \pi}$,
$$Q \models \exists {\vec {\bar \gamma}} < \kappa \ \kappa = 
\tau({\vec {\bar \gamma}},{\vec \epsilon}) {\rm , }$$
which contradicts $\kappa \notin H^Q_\Gamma(\kappa)$.
 
\hfill $\square$ (Claim \ref{claim1})

\bigskip
We may now define a strictly monotone $\varphi \colon \zeta+1 \rightarrow
\theta+1$ as follows. For $\alpha+1 \leq \zeta$, we let
$\varphi(\alpha)$ be the unique $\beta$ such that 
${\rm crit}(E_\beta^{\cal U}) =
\pi({\rm crit}(E_\alpha^{\cal T}))$. Moreover, we 
set $\varphi(\zeta) = \theta$.

We now wish to define, for $\beta \leq \theta$, 
elementary embeddings $\pi_\beta 
\colon {\cal M}_{\varphi^-(\beta)}^{\cal T} \rightarrow {\cal
M}_\beta^{\cal U}$
by  
$$\pi_\beta = (\pi_{\beta \theta}^{\cal U})^{-1} \circ \pi \circ
\pi_{\varphi^-(\beta) \zeta}^{\cal T}.$$
In order to see that these maps are well-defined it
suffices to prove the following.

\begin{claim}\label{claim2}
Let $\beta \leq \theta$, and set $\alpha = \varphi^-(\beta)$. Then 
${\rm ran}(\pi \circ
\pi_{\alpha \zeta}^{\cal T}) \subset {\rm ran}(\pi_{\beta \theta}^{\cal U})$.
\end{claim}

{\sc Proof} of Claim \ref{claim2}. 
Let $x \in {\cal M}_\alpha^{\cal T}$. Let $\nu = \mu^{(+3)W}$ for some $\mu$ 
such that $\pi_{\alpha \zeta}^{\cal T}(x) \in W||\nu$.
Let $Q$ and $\Gamma$ be such that
$Q \triangleright W||\nu$ is a universal weasel and $\Gamma$ is thick in
$Q$.
Let $E = E_{\pi \upharpoonright (W||\nu)}$ be the long extender derived
from $\pi \upharpoonright (W||\nu)$, and let
$${\tilde \pi} \colon Q \rightarrow_E {\tilde Q} = Ult(Q;E)$$
be the ultrapower of $Q$ by $E$. 
Again, $\pi \upharpoonright (W||\nu) = {\tilde \pi} 
\upharpoonright (W||\nu)$, ${\tilde Q}$ is universal,
${\tilde \pi}(\nu) = \pi(\nu)$ by Lemma \ref{folklore2}, 
and $W'||\pi(\nu) \triangleleft {\tilde Q}$.
Let $${\tilde \Gamma} = \Gamma \cap \{ \epsilon \colon {\tilde \pi}(\epsilon) =
\epsilon \}.$$ We'll have that ${\tilde \Gamma}$ is
thick in both $Q$ and ${\tilde Q}$.

Now let $\pi_{\alpha \zeta}^{\cal T}(x) = \tau^{Q}({\vec \gamma},{\vec
\epsilon})$, where $\tau$ is a Skolem term,
$\pi_{\alpha \zeta}^{\cal T} \upharpoonright {\vec \gamma} = {\rm id}$, and ${\vec
\epsilon} \in {\tilde \Gamma}$. 

Set $y = \pi \circ \pi_{\alpha \zeta}^{\cal T}(x)$.
Notice that $y \in W'||\pi(\nu)$. Let $\beta'$ be least with $\pi_{\beta' \theta}^{\cal
U} \upharpoonright \pi(\nu) = {\rm id}$. There is then an elementary embedding
$$\sigma \colon {\cal M}_{\beta'}^{\cal U} \rightarrow {\tilde Q}$$ with $\sigma
\upharpoonright \pi(\nu) = {\rm id}$.
We now get that
$y = 
{\tilde \pi} \circ \pi_{\alpha \zeta}^{\cal T}(x)
= \tau^{{\tilde Q}}(\pi({\vec \gamma}),{\vec
\epsilon})$, where $\pi_{\beta \theta}^{\cal U} \upharpoonright 
\pi({\vec \gamma}) = {\rm id}$. In particular,
$y \in \sigma \circ {\rm ran}(\pi_{\beta \beta'}^{\cal U})$. But $\sigma
\upharpoonright {\rm TC}(\{ y \}) = id$, and hence
$y \in {\rm ran}(\pi_{\beta \beta'}^{\cal U})$. But 
$\pi_{\beta' \theta}^{\cal U} \upharpoonright {\rm TC}(\{ y
\}) = {\rm id}$, and therefore 
$y \in {\rm ran}(\pi_{\beta \theta}^{\cal U})$.

\hfill $\square$ (Claim \ref{claim2})

\bigskip
We are now left with having to verify that (a) through (d) as in the statement of
Theorem \ref{mainthm} hold. Note that (a) and (b) are trivial. Let us show (c) and (d).
It is easy to see that the following suffices to establish (c).

\begin{claim}\label{claim3}
Let $\alpha+1 \leq \zeta$, and set $\kappa = {\rm crit}(E_\alpha^{\cal T})$. 
Then $\pi_{\varphi(\alpha)+1}(\kappa) = \pi_{\varphi(\alpha)}(\kappa)$. 
\end{claim}

{\sc Proof} of Claim \ref{claim3}.
Let us write $\beta = \varphi(\alpha)$. We know that 
$\pi_{\beta+1}(\kappa) = \pi(\kappa) = {\rm crit}(E_\beta^{\cal U})$.
Setting $\lambda = \pi(\kappa)$, we must therefore prove that
\begin{eqnarray}\label{lemma3-1}
\pi_{\beta \theta}^{\cal U}(\lambda) = \pi(\pi_{\alpha \zeta}^{\cal T}(\kappa)).
\end{eqnarray}

Let $\alpha' \in {\rm OR}$ be such that $\pi_{\alpha' \zeta}^{\cal T} \upharpoonright
(\pi_{\alpha \alpha'}^{\cal T}(\kappa)+1) = {\rm id}$.
Let $\beta' \in {\rm OR}$ be such that $\pi_{\beta' \theta}^{\cal U} \upharpoonright
(\pi_{\beta \beta'}^{\cal U}(\lambda)+1) = {\rm id}$
and $\pi_{\beta' \theta}^{\cal U} \upharpoonright (\pi(\pi_{\alpha \alpha'}^{\cal
T}(\kappa))+1) = {\rm id}$.
(If $\zeta \in {\rm OR}$ rather than $\zeta = {\rm OR}$ then
we may just let $\alpha' = \zeta$;
similarily, if $\theta \in {\rm OR}$ rather than $\theta = {\rm OR}$ then
we may just let $\beta' = \theta$.)
Let
$\nu = \mu^{(+3)K}$ for some $\mu$ 
be such that 
${\cal T} \upharpoonright (\alpha'+1)$ as well as 
${\cal U} \upharpoonright (\beta'+1)$
both
``live on'' $K||\nu$, i.e., such
that ${\rm lh}(E_\gamma^{\cal T}) < \pi_{0 \gamma}^{\cal T}(\nu)$ whenever $\gamma+1
\leq \alpha'$
and ${\rm lh}(E_\gamma^{\cal U}) < \pi_{0 \gamma}^{\cal U}(\nu)$ whenever $\gamma+1
\leq \beta'$. We have that
$$\pi_{0 \alpha'}^{\cal T}(\nu) > 
\pi_{\alpha \alpha'}^{\cal T}(\kappa) = \pi_{\alpha \zeta}^{\cal T}(\kappa) {\rm \ \ \
and }$$
$$\pi_{0 \beta'}^{\cal U}(\nu) > 
\pi_{\beta \beta'}^{\cal U}(\lambda) = \pi_{\beta \lambda}^{\cal U}(\lambda).$$
We may and shall in fact assume that 
$$\pi_{\alpha' \zeta}^{\cal T} \upharpoonright \pi_{0 \alpha'}^{\cal T}(\nu) = {\rm
id} {\rm \ \ \ and }$$
$$\pi_{\beta' \theta}^{\cal U} \upharpoonright \pi_{0 \beta'}^{\cal U}(\nu) = {\rm
id}.$$

Let $Q \triangleright K||\nu$ be a very soundness witness for $K||\nu$. 
We may construe 
${\cal T} \upharpoonright (\alpha'+1)$ and ${\cal U} \upharpoonright (\beta'+1)$
as iteration trees
acting on $Q$ rather than on $K$. More precisely,
we let ${\cal T}^*$ be the iteration of $Q$ which is such that
${\rm lh}({\cal T}^*) = \alpha'+1$ and
$E_\gamma^{{\cal T}^*} = E_\gamma^{\cal T}$ whenever $\gamma < \alpha'+1$,
and we let
${\cal U}^*$ be the iteration of $Q$ which is such that
${\rm lh}({\cal U}^*) = \beta'+1$ and
$E_\gamma^{{\cal U}^*} = E_\gamma^{\cal U}$ whenever $\gamma < \beta'+1$.
We let $${\tilde \pi} \colon {\cal M}_{\alpha'}^{{\cal T}^*} \rightarrow 
{\cal R} = Ult({\cal M}_{\alpha'}^{{\cal T}^*};\pi \upharpoonright
(W||\pi_{0 \alpha'}(\nu)).$$ We'll have that
$\pi \upharpoonright
(W||\pi_{0 \alpha'}(\nu)) = {\tilde \pi} \upharpoonright
(W||\pi_{0 \alpha'}(\nu))$, 
${\cal R}$ is a universal weasel,
${\tilde \pi}(\pi_{0 \alpha'}(\nu)) =
\pi(\pi_{0 \alpha'}(\nu))$ by Lemma \ref{folklore2}, 
and $W'||\pi(\pi_{0 \alpha'}(\nu))
\triangleleft {\cal R}$.
Let ${\tilde {\cal R}}$ be the common co-iterate of ${\cal R}$ and ${\cal M}^{{\cal
U}^*}_{\beta'}$, and let $$i \colon {\cal R} \rightarrow {\tilde {\cal R}} {\rm \ \ \
and }$$ $$j \colon {\cal M}^{{\cal
U}^*}_{\beta'} \rightarrow {\tilde {\cal R}}$$ be the maps arising from the
coiteration.
Notice that $$i \upharpoonright (\pi(\pi_{\alpha \alpha'}^{\cal T}(\kappa)+1) = {\rm
id} {\rm \ \ \ and }$$
$$j \upharpoonright (\pi_{\beta \beta'}^{\cal U}(\lambda)+1) = {\rm id }.$$
Finally, let $\Gamma$ be thick in ${\cal M}_\alpha^{{\cal T}^*}$, 
${\cal M}_{\alpha'}^{{\cal T}^*}$, ${\cal R}$, ${\cal M}^{{\cal
U}^*}_\beta$, ${\cal M}^{{\cal
U}^*}_{{\beta}'}$, and ${\tilde {\cal R}}$ and such that for all
$\epsilon \in \Gamma$ we have that $$\pi_{\alpha \alpha'}^{{\cal
T}^*}(\epsilon) = {\tilde \pi}(\epsilon) = i(\epsilon) = \pi_{\beta \beta'}^{{\cal
U}^*}(\epsilon) = j(\epsilon) = \epsilon.$$

Now let $\kappa = \tau^{{\cal M}_\alpha^{{\cal T}^*}}({\vec \gamma},{\vec \epsilon})$, 
where
${\vec \gamma} < \kappa$ and ${\vec \epsilon} \in \Gamma$.
In order to show (\ref{lemma3-1}) it suffices to prove that
\begin{eqnarray}\label{lemma3-2}
\lambda = \tau^{{\cal M}_\beta^{{\cal U}^*}}(\pi({\vec \gamma}),{\vec \epsilon}) 
{\rm , }
\end{eqnarray}
because then $\pi(\pi_{\alpha \zeta}^{\cal T}(\kappa)) = 
\pi(\pi_{\alpha \alpha'}^{\cal T}(\kappa)) = 
{\tilde \pi}(\pi_{\alpha \alpha'}^{{\cal T}^*}(\kappa)) =
i({\tilde \pi}(\pi_{\alpha \alpha'}^{{\cal T}^*}(\kappa))) =
\tau^{{\tilde {\cal R}}}(\pi({\vec \gamma}),{\vec \epsilon}) =
j(\pi_{\beta \beta'}^{{\cal U}^*}(\tau^{{\cal M}_\beta^{{\cal U}^*}}(\pi({\vec \gamma}),{\vec \epsilon})))
= j(\pi_{\beta \beta'}^{{\cal U}^*}(\lambda)) =
j(\pi_{\beta \beta'}^{\cal U}(\lambda)) =
\pi_{\beta \beta'}^{\cal U}(\lambda) = \pi_{\beta \theta}^{\cal U}(\lambda)$.

Write $\lambda = {\bar \tau}^{{\cal M}_\beta^{{\cal U}^*}}({\vec \gamma}_0,{\vec
\epsilon}_0)$, where ${\vec \gamma}_0 < \lambda$ and ${\vec \epsilon}_0 \in \Gamma$.

Suppose that $\lambda < 
\tau^{{\cal M}_\beta^{{\cal U}^*}}(\pi({\vec \gamma}),{\vec \epsilon}).$
Then, using 
$j \circ \pi_{\beta \beta'}^{{\cal U}^*}$, 
$\lambda \leq {\bar \tau}^{{\tilde {\cal R}}}({\vec \gamma}_0,{\vec
\epsilon}_0) < \tau^{\tilde {\cal R}}(\pi({\vec \gamma}),{\vec \epsilon})$, 
and therefore
\begin{eqnarray}\label{lemma3-3}
{\tilde {\cal R}} \models \exists {\vec {\bar \gamma}} < \lambda \ ( \lambda 
\leq {\bar \tau}({\vec {\bar \gamma}},{\vec
\epsilon}_0) < \tau(\pi({\vec \gamma}),{\vec \epsilon})).
\end{eqnarray}
Using $i \circ {\tilde \pi}$, hence,
\begin{eqnarray}\label{lemma3-4}
{\cal M}_{\alpha'}^{{\cal T}^*} \models \exists {\vec {\bar \gamma}} < \kappa \ ( \kappa 
\leq {\bar \tau}({\vec {\bar \gamma}},{\vec
\epsilon}_0) < \tau({\vec \gamma},{\vec \epsilon})).
\end{eqnarray}
But (\ref{lemma3-4}) is wrong, as 
$\tau^{{\cal M}_{\alpha'}^{{\cal T}^*}}({\vec \gamma},{\vec \epsilon}) =
\pi_{\alpha \zeta}^{\cal T}(\kappa)$ is the least ordinal 
in $H^{{\cal M}_{\alpha'}^{{\cal T}^*}}(\kappa) \setminus \kappa$. Therefore, 
$\lambda \geq 
\tau^{{\cal M}_\beta^{{\cal U}^*}}(\pi({\vec \gamma}),{\vec \epsilon}).$

However, if $\lambda > 
\tau^{{\cal M}_\beta^{{\cal U}^*}}(\pi({\vec \gamma}),{\vec \epsilon})$
then $\tau^{\tilde {\cal R}}(\pi({\vec \gamma}),{\vec \epsilon}) =
\tau^{{\cal M}_\beta^{{\cal U}^*}}(\pi({\vec \gamma}),{\vec \epsilon}) < \lambda$, and thus
by using $i \circ {\tilde \pi}$,
$\kappa = \tau^{{\cal M}_\alpha^{{\cal T}^*}}({\vec \gamma},{\vec \epsilon}) < \kappa$, which is nonsense.
We have therefore established that (\ref{lemma3-2}) holds.

\hfill $\square$ (Claim \ref{claim3})

\begin{claim}\label{claim4} Let $\alpha+1 \leq \zeta$, and set $\beta = \varphi(\alpha)$.
Then $\pi_\beta {\rm " } E_\alpha^{\cal T} \subset E_\beta^{\cal U}$.
In fact, $\pi_\beta(E_\alpha^{\cal T})= E_\beta^{\cal U}$.
\end{claim}

{\sc Proof} of Claim \ref{claim4}.
Set $\kappa = {\rm crit}(E_\alpha^{\cal T})$, and
$\lambda = {\rm crit}(E_\beta^{\cal U}) = \pi(\kappa)$. 
Fix $X \in E_\alpha^{\cal T}$. We then get $\kappa \in \pi_{\alpha
\zeta}^{\cal T}(X)$, hence $\lambda = \pi(\kappa) \in \pi(\pi_{\alpha
\zeta}^{\cal T}(X)) = \pi_{\beta \theta}^{\cal U}(\pi_\beta(X))$.
But this means that $\pi_\beta(X)
\in E_\beta^{\cal U}$.
This shows that $\pi_\beta {\rm " } 
E_\alpha^{\cal T} \subset E_\beta^{\cal U}$.

However, we may now apply Lemma \ref{observation} and deduce that
in fact $\pi_\beta(E_\alpha^{\cal T})= E_\beta^{\cal U}$.

\hfill $\square$ (Lemma \ref{claim4})

This finishes the proof of Theorem \ref{main-theorem}.

\hfill $\square$ (Theorem \ref{main-theorem})

\section{The main results.}

In order to obtain our main results we shall make use of 
the following criterion.

\begin{lemma}\label{criterion}
Let $A \subset {}^\omega \omega$. Then
$A$ is coanalytic if and only if there is a
function $n \colon {}^{< \omega} \omega \rightarrow \omega$
with $n(t) \geq n(s) \geq {\rm lh}(s)$ for $s \subset t \in {}^{<
\omega} \omega$ and a system $(\sigma_{s t} \colon s \subset t \in {}^{<
\omega} \omega)$ such that for $s \subset t \in {}^{<
\omega} \omega$, $\sigma_{s t} \colon n(s) \rightarrow n(t)$ is
order-preserving, and
for all $x \in {}^\omega \omega$,
$$x \in A \Longleftrightarrow {\rm lim \ dir} 
((n(s) \colon s \subset x),
(\sigma_{s t} \colon s \subset t \subset x))
{\rm \ is \ wellfounded.}$$
\end{lemma}

{\sc Proof} of Theorem \ref{thm1}. 
Let $A \subset {}^\omega \omega$, and let 
$$((M_s \colon s \in {}^{< \omega}\omega),(\pi_{st} 
\colon s \subset t \in {}^{< \omega}\omega))$$ be an 
$\omega$-closed embedding normal form for $A$.
We have that $M_0 = V = K$. 
For an arbitrary $s \in {}^{<\omega}\omega$,
$M_s$ must be a {\em finite}
iterate of $K$, because ${}^\omega M_s \subset M_s$.
Let 
${\cal T}_s = (E^{{\cal T}_s}_i \colon i<n(s))$ denote the iteration from 
$K$ to $M_s$. In particular,
$n(s)$ is the number of critical points of this
iteration from $K$ to $M_s$.
Set $\kappa_i^s = {\rm crit}(E^{{\cal T}_s}_i)$ for $i<n(s)$.
By Theorem \ref{main-theorem}, we know that in particular
$$\pi_{st} {\rm " } \{ \kappa_0^s, ..., \kappa_{n(s)-1}^s \}
\subset \{ \kappa_0^t, ..., \kappa_{n(t)-1}^t \}$$ whenever
$s \subset t \in {}^{<\omega} \omega$. Hence if
$s \subset t \in {}^{<\omega} \omega$ then $\pi_{s t}$ induces
an order-preserving map 
$\sigma_{s t}$ such that $\pi_{s t}(\kappa_i^s) = \kappa_{\sigma_{s t}(i)}^t$
for all $i < n(s)$. 

We claim that $((n(s) \colon s \in {}^{<\omega} \omega),
(\sigma_{s t} \colon s \subset t \in {}^{<\omega} \omega))$ witnesses
that $A$ is coanalytic. The non-trivial part here is to show that
if the direct limit of $((n(s) \colon s \subset x),
(\sigma_{s t} \colon s \subset t \subset x))$ is wellfounded
then the direct limit of $((M_s \colon s \subset x),
(\pi_{s t} \colon s \subset t \subset x))$ is also wellfounded.

Let 
$x \in {}^\omega \omega$ be such that 
the direct limit of $((n(s) \colon s \subset x),
(\sigma_{s t} \colon s \subset t \subset x))$ is wellfounded,
and let $\theta < \omega_1$ be its ordertype.
The point is that we may now use Theorem \ref{main-theorem}
to construct an iteration ${\cal T}$ of $K$ 
of length $\theta+1$ together with commuting elementary embeddings
from the models ${\cal M}^{{\cal T}_{x \upharpoonright n}}_i$ into the
models ${\cal M}^{\cal T}_j$. 
This construction will in particular give
embeddings $\pi_{s x} \colon M_s \rightarrow {\cal M}^{\cal T}_\theta$
such that $\pi_{t x} \circ \pi_{s t} = \pi_{s x}$ whenever
$s \subset t \subset x$. We leave the straightforward details of this
construction to the 
reader. But then the direct limit of $((M_s \colon s \subset x),
(\pi_{s t} \colon s \subset t \subset x))$
can be embedded into (in fact, is equal to!)
${\cal M}^{\cal T}_\theta$, so that it must be wellfounded.

\hfill $\square$ (Theorem \ref{thm1}) 

\begin{lemma}\label{iterations-are-short}
{\rm ($\lnot \ 0^{\rm long}$)} Let $\pi \colon V \rightarrow M$, where
$M$ is transitive and ${}^\omega M \subset M$. Then $K^M$ is a finite iterate
of $K$.
\end{lemma}

{\sc Proof.} By \cite[Theorem 3.23]{peter} (cf.~also \cite[Theorem 1.3]{bill}
and the remark right after it), there is a maximal 
Prikry system ${\mathbb C}$ for $K$. 
Let us construe ${\mathbb C}$ as a set of ordinals; we shall have
that if
$\nu$ denotes the order type of the measurable cardinals 
of $K$ then ${\mathbb C}$ will be
of order type at most $\omega \cdot \nu$.
By $\lnot \ 0^{\rm long}$, $\nu$ is less than the least measurable cardinal
of $K$, which in turn is less than or equal to ${\rm crit}(\pi)$.
The covering lemma (cf.~\cite[Theorem 3.23]{peter})
says that for each set $X$ of ordinals
there is some function $f \in K$ and some $\alpha < \aleph_2
\cdot {\rm Card}(X)^+$ such that $X \subset f {\rm " } (\alpha \cup
{\mathbb C})$.

By the Dodd-Jensen Lemma, $K^M$ is universal and hence a non-dropping 
iterate of $K$.
Now suppose that $K^M$ would be an infinite iterate of $K$. 
Let $A$ be the set of the first $\omega$ many critical points 
of the iteration from $K$ to $K^M$. By elementarity,
there is some function $f \in K^M$ and some $\alpha < \aleph_2$ 
such that $A \subset f {\rm " } (\alpha \cup
\pi({\mathbb C}))$. However, $\pi({\mathbb C}) =
\pi {\rm " }{\mathbb C}$, as ${\rm otp}({\mathbb C}) < {\rm crit}(\pi)$.
That is, $A \subset f {\rm " } (\alpha \cup
\pi {\rm " }{\mathbb C})$. 

Let $f = \pi({\bar f})(a)$, where ${\bar f} \in K$ and $a$ is a finite set
of critical points of the iterations from $K$ to $K^M$.
Let $\kappa \in A \setminus a$. Then $\kappa = (\pi({\bar f})(a))(\xi)$,
where $\xi \in \alpha \cup \pi {\rm " }{\mathbb C}$, so that
in particular $\kappa$
is in the hull of ${\rm ran}(\pi) \cup a$ formed inside $K^M$.
But this is a contradiction, as $\kappa$ is one of the critical points
of the iteration from $K$ to $K^M$.

\hfill $\square$ (Lemma \ref{iterations-are-short})

\bigskip
{\sc Proof} of Theorem \ref{thm2}. This can now be shown exactly as
Theorem \ref{thm1} above, using Lemma \ref{iterations-are-short}.

\hfill $\square$ (Theorem \ref{thm2})

\end{document}